\documentclass[a4paper
	       ]{scrartcl} 
\usepackage{typearea}
\typearea[current]{calc}  

\usepackage[latin1]{inputenc}
\usepackage[ngerman, english]{babel}

\usepackage{amssymb}
\usepackage{amsmath}
\usepackage{amsthm}
\usepackage{mathtools}

\usepackage{graphicx} 
\usepackage{layout}

\title{The Pólya branching process and limit theorems for conditioned random fields}
\author{Benjamin Nehring\thanks{Universit\"at Potsdam, Institut f\"ur Mathematik, Am Neuen Palais 10, D-14469 Potsdam, Germany; e-mail {\sf benjamin.nehring@gmail.com}} \and Mathias Rafler\thanks{TU M\"unchen, Zentrum Mathematik M5, Boltzmannstr. 3, D-85747 Garching bei M\"unchen, Germany; e-mail {\sf rafler@ma.tum.de}}\\
}

\renewcommand{\Pr}{\mathbb P} 

\newcommand{\Poi}{\mathbf P} 

\newcommand{\Poy}{\mathsf S\!} 
\newcommand{\Qp}{\mathsf Q} 
\newcommand{\PoyS}[2][\ast]{\mathsf S^{#1}_{\! #2}} 

\newcommand{\Pp}{\mathsf P} 
\newcommand{\GP}{\mathsf D} 

\newcommand{\B}{\mathcal{B}} 
\newcommand{\Bbd}{\mathcal{B}_{0}} 

\newcommand{\Fsig}{\mathcal{F}}
\newcommand{\Esig}{\mathcal{E}}

\def\Ftail{\Fsig_{\infty}}
\def\Etail{\Esig_{\infty}}

\def\Ef{\mathbb E}

\def\fsm{\text{-a.s.}}

\newcommand{\M}{\mathcal M}
\newcommand{\MX}{\mathcal M(X)}

\newcommand{\Mpm}{\mathcal M^{\cdot\cdot}}
\newcommand{\MpmX}{\mathcal M^{\cdot\cdot}(X)}

\newcommand{\N}{\mathbb N}

\newcommand{\R}{\mathbb R}

\renewcommand{\d}{\mathrm{d}}

\newcommand{\supp}{\operatorname{supp}}

\renewcommand{\exp}{\operatorname{exp}}
\newcommand{\e}{\operatorname{e}}

\newcommand{\bigast}{\mathop{\scalebox{1.8}{\raisebox{-0.2ex}{$\ast$}}}}%

\renewcommand{\phi}{\varphi}
\renewcommand{\theta}{\vartheta}

\newcommand{\F}{\mathcal{F}}

\newcommand{\eqa}[1]{\begin{align*}#1\end{align*}}
\newcommand{\equ}[1]{\begin{equation*}#1\end{equation*}}

\newcommand{\equn}[1]{\begin{equation}#1\end{equation}}

\theoremstyle{definition}
\newtheorem{defdefinition}{Definition}[section]

\theoremstyle{plain}
\newtheorem{defsatz}[defdefinition]{Theorem}
\newtheorem{defsatzdef}[defdefinition]{Theorem and Definition}
\newtheorem{defprop}[defdefinition]{Proposition}
\newtheorem{deflemma}[defdefinition]{Lemma}
\newtheorem{deffolgerung}[defdefinition]{Corollary}
\newtheorem{defbeispiel}[defdefinition]{Example}

\theoremstyle{remark}
\newtheorem{defbemerkung}[defdefinition]{Remark}

\newcommand{\defn}[2]{\begin{defdefinition}[#1]#2\end{defdefinition}}

\newcommand{\satz}[1]{\begin{defsatz}#1\end{defsatz}}
\newcommand{\satzn}[2]{\begin{defsatz}[#1]#2\end{defsatz}}

\newcommand{\prop}[1]{\begin{defprop}#1\end{defprop}}

\newcommand{\bem}[1]{\begin{defbemerkung}#1\end{defbemerkung}}
\newcommand{\lemma}[1]{\begin{deflemma}#1\end{deflemma}}

\newcommand{\korollar}[1]{\begin{deffolgerung}#1\end{deffolgerung}}

\newcommand{\beispiel}[1]{\begin{defbeispiel}#1\end{defbeispiel}}

\newcounter{margcount}
\numberwithin{equation}{section}

\begin{document}

\maketitle

\begin{abstract}
The first aim is to construct generalizations of Pólya type point process by applying a branching mechanism to these point processes. Conditions are given under which these point processes satisfy an integration by parts formula. Furthermore we compute their Palm kernels, which turn out to be superpositions of different point processes. Secondly we identify all point processes whose local characteristics agree with these of a fixed branching of a Pólya type point process as mixtures of branchings of Pólya type point process and show that in this case also they are characterized by an integration by parts formula.
\end{abstract}

%

\section{Introduction}

The fundamental example of a point process being a solution of a partial integration formula is the Poisson process, which is characterized by Mecke's celebrated formula
\begin{equation} \label{eq:intro:mecke}
  C_{\Pr}(h)=\iint h(x,\mu+\delta_x)\rho(\d x)\Pr(\d\mu)
\end{equation}
saying firstly that the Poisson process with intensity measure $\rho$ satisfies equation~\eqref{eq:intro:mecke}, and secondly that this functional equation has a unique solution, the Poisson process with intensity measure $\rho$.

Point processes satisfying such a disintegration formula have been considered e.g. in~\cite{MWM79,NZ79,oK78}, giving conditions under which a disintegration of the Campbell measure is available. In general $\rho$ has to be replaced by some kernel depending on $\mu$. Starting from Mecke's formula one might try to generalize the Poisson process by replacing the $\rho$ by some simple kernel questioning the existence of such a process. Indeed, in~\cite{hZ09} Zessin started with the kernel $\pi(\mu,\,\cdot\,)=z(\rho+\mu)$, which implements a reinforcement scheme, and showed that there exists a unique point process which satisfies the modified disintegration formula. This point process he called the Pólya sum process. For the moment $z$ is a positive number stricly smaller than 1, an interpretation is given later.

A natural generalization of the Poisson process is the Cox process, where the intensity measure is considered as a random measure. Via this generalization, the Poisson process is directly connected with the Pólya sum process, when viewing from a Bayesian standpoint, see~\cite{GW82,mR11as}: If the directing measure, the prior, is a Poisson gamma random measure, then the posterior is again a Poisson gamma random measure and the Bayes estimator of the directing intensity measure is exactly the kernel of the Pólya sum process. Thus by using techniques from Bayesian statistics, the Pólya sum process turns out to be a particular Cox process.

Considering the construction of the Pólya sum process, it turns out that this point process directly generalizes Hoppe's urn scheme~\cite{fH84}: Given a realization $\mu$ of points, the intensity for a new point is $\pi(\mu,\,\cdot\,)=z(\rho+\mu)$, the Pólya sum kernel, thus points in a configuration always get a reward. We want to discuss several possibilities to generalize the Pólya sum process starting from a modification of the kernel $\pi$. First of all, one might replace the "+" by a "-", which yields the Pólya difference process as considered in~\cite{NZ11}, or secondly the fixed number $z$ may be replaced by some measurable function. Indeed, both situations occur for finite systems of quantum particles. In~\cite{BZ11}, bosonic particle systems are related to the Pólya sum process and fermionic particles to the Pólya difference process and that takes the role of the distribution of a single particle.


While these generalizations do not need further construction techniques apart from~\cite{hZ09}, it is more difficult to find an answer when replacing $\mu$ by a kernel $V_\mu$. Instead of rewarding the exact locations of the points, the idea is to reward regions which are in a to be specified sense related to the positions of the points. The aim is to obtain a disintegration formula of kind~\eqref{eq:intro:mecke} for the resulting point process, i.e. to obtain a Papangelou process. It turns out, that this process can be constructed from the Pólya sum process by a random motion of the points of the configuration of the Pólya sum process. We show that $V$ has to be constructed from a conditional probability. The result will be an infinitely divisible point process with a weak indepence of increments property. This process we call \emph{branching Pólya sum process}. An analogue result holds when starting with the Pólya difference process. Branching mechanism and the construction of the branching processes are subject of sections~\ref{sect:branchmech:branch}~and~\ref{sect:construction}, the technical parts and the proofs are postponed to section~\ref{sect:constrproofs}. In section~\ref{sect:properties} we compute the Palm kernels of the branching Pólya sum process, which yields an interpretation for $z$.

In a second step we characterize those point processes which are locally distributed like a conditioned branching Pólya sum process. The technique goes back to Dynkin~\cite{eD71a,eD71b}, see also~\cite{eD78,hF75}, and was used in~\cite{hZ76} to characterize mixed Poisson processes as canonical Gibbs states of a certain local specification, later in~\cite{mR11jtp} to identify the families of mixed Pólya sum processes as Gibbs states for particular local specifications. We show that in the situation of our conditioning we get a certain family of mixed branching Pólya processes. Moreover, from this analysis we obtain that these mixed branching processes are Papangelou processes themselves, but not characterized by this integration by parts formula in contrast to the branching Pólya processes. These aspects are part of section~\ref{sect:limits}, their proofs are contained in section~\ref{sect:limitsproofs}. As the Martin-Dynkin boundary has not been discussed so far, we added this discussion in section~\ref{sect:pdp}.

\section{The branching mechanism\label{sect:branchmech}}

Let $X$ be a Polish space with its Borel sets $\B$ and its bounded Borel sets $\B_0$. Denote by $\MX$ and $\MpmX$ the set of locally finite measures and point measures on $X$, respectively. When equipped with the $\sigma$-field generated by the evaluation mappings, or counting mappings in case of point measures, $\zeta_B(\mu)=\mu(B)$, then $\MX$ and $\MpmX$ are Polish spaces themselves. A random measure is a probability measure on $\MX$ and a point process if is is concentrated on $\MpmX$. In general we do not deal with simple point pocesses.

For a measurable, non-negative function $f$ write $\mu(f)$ for the integral of $f$ with respect to $\mu$. Let $\kappa$ be a kernel from $X$ to $X$, partially we denote by $\kappa(f)$ the mapping $x\mapsto \kappa_x(f)=\int f(y)\kappa_x(\d y)$ and for some $\rho\in\MX$ by $\kappa\rho$ the measure $\int \kappa_x\,\rho(\d x)$.

\subsection{Papangelou processes\label{sect:branchmech:papa}}

A Papangelou process $\Pp$ is a point process whose Campbell measure $C_\Pp$ on the lhs. satisfies the integration by parts formula
\begin{equation} \label{eq:branchmech:papa:ibpf}
  \iint h(x,\mu)\mu(\d x)\Pp(\d\mu)=\iint h(x,\mu+\delta_x)\pi(\mu,\d x)\Pp(\d\mu)
\end{equation}
for all measurable functions $h:X\times\MpmX\to[0,+\infty]$ for some kernel $\pi$ from $\MpmX$ to $X$. $\pi$ is called \emph{Papangelou kernel} and $\pi(\mu,B)$ can be understood as the expected number of points in $B\in\B$ given a realized point configuration $\mu$. 

Conditions on a point process $\Pp$ to satisfy an integration by parts formula are given e.g. in~\cite{MWM79,NZ79}. In this case, $\pi$ is unique, and moreover the iterated Papangelou kernels
\begin{equation*}
  \pi^{(m)}(\mu;\d x_1,\ldots,\d x_m)=\pi(\mu+\delta_{x_1}+\ldots+\delta_{x_{m-1}},\d x_m)\pi^{(m-1)}(\mu;\d x_1,\ldots,\d x_{m-1}),
\end{equation*}
$\pi^{(1)}=\pi$, are symmetric measures. It is immediate to see that this symmetry is equivalent to the \emph{cocycle condition}
\begin{equation}
  \pi(\mu+\delta_x,\d y)\pi(\mu,\d x)=\pi(\mu+\delta_y,\d x)\pi(\mu,\d y),
\end{equation}
i.e. to the symmetry of $\pi^{(2)}$.

Starting with a kernel $\pi$, the functional equation~\eqref{eq:branchmech:papa:ibpf} may not uniquely define a point process. However, Poisson process, Pólya sum process and Pólya difference process are characterized by $\pi(\mu,\,\cdot\,)=\rho$, $\rho\in\MX$, $\pi(\mu,\,\cdot\,)=z(\rho+\mu)$, $\rho\in\MX$ and $z\in(0,1)$, and $\pi(\mu,\,\cdot\,)=z(\rho-\mu)$, $\rho\in\MpmX$ and $z>0$. As a by-product we construct kernels for which equation~\eqref{eq:branchmech:papa:ibpf} has infinitely many solutions.

\subsection{Branching\label{sect:branchmech:branch}}

Let $\kappa$ be a kernel from $X$ to $X$. As a stochastic kernel $\kappa$ will take the role of a random motion of points of a point configuration $\mu\in\MpmX$. If we let $\Delta_y=\delta_{\delta_y}$ be the point process that realizes a single point at $y\in X$, then for $\mu\in\MpmX$ define
\equn{ \label{eq:setup:V-def}
  V_\mu=\underset{x\in \mu}{\bigast}\int \Delta_y\, \kappa_x(\d y),
}%
i.e. each point $x\in\mu$ is moved independently of all other points according to $\kappa_x$.

\defn{Branching processes}{
For $\mu\in\MpmX$, if $V_\mu$ given by equation~\eqref{eq:setup:V-def} is a locally finite point measure, $V_\mu$ is the \emph{branching} of the configuration $\mu$. Its mixture $\Poi$ with respect to a point process $\Pp$ is the \emph{branching} of $\Pp$,
\equ{
  \Poi=\int V_\mu\,\Pp(\d\mu).
}
If $\Pp=\Poy_{z,\rho}$ is the Pólya sum process, its branching is the \emph{branching Pólya sum process}, and if $\Pp=\GP_{z,\rho}$ is the Pólya difference process, its branching is the \emph{branching Pólya difference process}.} 

Successively we put assumptions on $\pi$ to firstly construct branching Pólya sum and branching Pólya difference processes, and secondly to discuss canonical Gibbs states in this context.
\begin{enumerate}
  \item[${\mathbf A}_1$] $\kappa$ is a stochastic kernel.
\end{enumerate}
Under ${\mathbf A}_1$, if $\mu=\delta_{x_1}+\ldots+\delta_{x_n}$ is a finite configuration of points and if $f$ is a non-negative, measurable function, then
\equ{
  V_\mu(f)=\int_{X^n}f(\delta_{y_1}+\ldots+\delta_{y_n})\kappa_{x_1}(\d y_1)\cdots\kappa_{x_n}(\d y_n).
}
Particularly points with multiplicity are considered as multiple points. Since in $\mu$ may be a locally finite but infinite configuration of points, $\kappa$ shall not be allowed to put an infinite number of points into any bounded region.

The branching $\Poi$ is a doubly stochastic point process which takes the realization of a directing point process and puts each of its points independently at different locations. Our focus lies on the underlying process being the Pólya sum process $\Poy_{z,\rho}$ and the Pólya difference process, but they may be replaced by the Poisson process or any other point process. In fact, as soon as the directing point process has independent increments, or some slightly weaker independence condition, most of the following analysis can be carried out.

For the moment assume that $\kappa$ is given in the following way. We will see and motivate in Lemma~\ref{thm:equiv-cocycle}, whose proof is postponed to section~\ref{sect:constrproofs}, this choice: Let $H$ be a probability measure on $(X,\B)$, $\Esig\subseteq\B$ be a sub-$\sigma$-algebra of $\B$ and $\kappa$ be a regular version of $H$ conditioned on $\Esig$,
\equn{ \label{eq:setup:kappa-def}
  \kappa_x(A)=H(A|\Esig)(x).
}
This choice has a strong influence on how point configurations are transformed, but at first consider some examples.
\beispiel{ \label{bsp:setup}
\begin{enumerate}
  \item Let $\Esig=\B$, then $\kappa_x=\delta_x$ is a regular version of $H(\,\cdot\,|\Esig)$.
  \item For $X=\R^d$ may $\Esig$ be set set of permutation invariant sets, the $\sigma$-field of \emph{symmetric events}.
  \item Also for $X=\R^d$, $\Esig$ may be the set of \emph{isotropic events}, i.e. $A\in\Esig$ if $A=OA$ for any orthogonal transformation $O$.
  \item Let $X_1,X_2,\ldots\in\Bbd$ be a locally finite partition of $X$ into bounded sets and $\Esig=\sigma(\{X_i:i\geq 1\})$. \label{bsp:setup:part}
  \item For $\Esig=\{\emptyset,X\}$ we get $\kappa_x(A)=H(A)$.
\end{enumerate}
Thus in the first case, points stay at their location and the branching process is the process itself. In the second case, a point $x$ is transformed by a permutation of its coordinates, and in the third case by a rotation around the origin. In the pre-last case $X$ is partitioned into separate islands and points are not allowed to migrate between them. We exclude the last example soon.
}

\bem{
Let $B\in\Esig$, then for $H$-a.e. $x\in B$ we have $\kappa_x(B)=1$ since $\kappa_x$ is given by a conditional probability, hence the image of $x$ is again contained in $B$. Thus the richer $\Esig$ is, the more restrictions are put on the branching mechanism. }

\lemma{ \label{thm:equiv-cocycle}
Let $\kappa$ be a stochastic kernel from $X$ to $X$ and define the \emph{branching} $\pi_\kappa$ of a kernel from $\M$ to $X$ by
\begin{equation*}
  \pi_{\kappa}(\mu,f)=\pi(\mu,\kappa(f))
\end{equation*} 
for non-negative, measurable $f$. If $\pi$ is the Pólya sum kernel $\pi(\mu,\d x)=z(\rho+\mu)(\d x)$ or the Pólya difference kernel $\pi(\mu,\d x)=z(\rho-\mu)(\d x)$, then the following statements are equivalent
\begin{enumerate}
  \item $\pi_{\kappa}$ satisfies the cocycle condition 
  \item $\kappa$ is a regular conditional probability.
\end{enumerate}
}%
The proof is postponed to section~\ref{sect:constrproofs} and we strengthen ${\mathbf A}_1$ to
\begin{enumerate}
  \item[${\mathbf A}_1'$] $\kappa$ is a regular version of a conditional probability.
\end{enumerate}

\subsection{The construction of the branching Pólya sum process \label{sect:construction}}

Let $\Esig_0=\B_0\cap\Esig$ be the set of bounded invariant events. From now on we ensure that $\Esig$ and $\Esig_0$ are rich enough:
\begin{enumerate}
  \item[${\mathbf A}_2$] There exists an increasing sequence $B_1\subset B_2\subset\ldots$ of sets in $\Esig_0$ such that $\bigcup_n B_n=X$ and for each $C\in\Bbd$ there exists $n\in\N$ sucht that $C\subseteq B_n$.
\end{enumerate}
This assumption ensures the existence of a locally finite partition $(X_{j})_{j\in\N}$ of $X$ into invariant bounded sets. Assumptions ${\mathbf A}_1'$ and ${\mathbf A}_2$ are sufficient for a branching of a Pólya process being a Papangelou process.
\satz{ \label{thm:papangelou}
Let $\pi$ be the Pólya sum or the Pólya difference kernel and $\kappa$ be a branching mechanism satisfying ${\mathbf A}_1'$, and assume ${\mathbf A}_2$. Then the branching $\Poi$ of the given Pólya process is the unique Papangelou process with Papangelou kernel $\pi_{\kappa}$, i.e. $\Poi$ satisfies for all non-negative, measurable functions $h$
\begin{equation} \label{eq:sect:construction:ibpf}
 C_{\Poi}(h)=\iint h(x,\mu+\delta_{x})\, \pi_{\kappa}(\mu,\d x) \Poi(\d\mu).
\end{equation}
}
Essentially we prove that the branching can be written as the superposition of independent, finite Papangelou processes. Thereby the main problem is to show that $\Poi$ solves equation~\eqref{eq:sect:construction:ibpf}, the existence of the point process is guaranteed by the clustering theorem 4.2.3 in~\cite{MKM78}: A cluster field on $X$ is a kernel $\Pi$ from $X$ to the set of point processes. Then
\begin{equation}\label{eq:constr:cluster}
  \Poi=\int\underset{x\in \mu}{\bigast} \Pi_{x}\,\Pp(\d\mu)
\end{equation}
exists as a point process if and only if its intensity measure is locally finite, i.e. for all $C\in\Bbd$
\begin{equation}\label{eq:constr:cond}
  \nu^{1}_{\Poi}(C)=\nu^{1}_{\Pp}(\nu^{1}_{\Pi_{(\cdot)}}(C))<\infty.
\end{equation}

In our situation, where $\Pp=\Poy_{z,\rho}$ or $\Pp=\GP_{z,\rho}$ and $\Pi_{x}=\int \Delta_{y}\, \kappa_{x}(dy)$, condition~\eqref{eq:constr:cond} turns into $\kappa\rho(C)<\infty$ for all $C\in\Bbd$. Thus if $(X_{j})_{j\in\N}$ is a locally finite partition of $X$ into bounded, invariant sets, then for every $C\in\Bbd$ only for finitely many $j\in\N$, $C\cap X_{j}\neq\emptyset$ and
\begin{equation*}
  \kappa\rho(C)=\sum_{j\in\N} \rho(1_{X_j}\kappa(C))=\sum_{j\in\N} \rho(\kappa(C \cap X_{j}))<\infty.
\end{equation*}
Hence under ${\mathbf A}_2$ the clustering~\eqref{eq:constr:cluster} exists.

What remains is to show that $\Poi$ is a Papangelou process with kernel $\pi_\kappa$. The basic step is contained in the following proposition:

\prop{ \label{thm:superpos}
Let $\Qp_j$ be the point process on $X$ be given by 
\begin{equation*}
  \Qp_j(\phi)=\frac{1}{\Xi_j} \sum_{n=0}^{\infty} \frac{1}{n!} \int_{X_j^n} \phi(\delta_{x_{1}}+\ldots+\delta_{x_{n}}) \, \pi_{\kappa}^{(n)}(0;\d x_{1},\ldots,\d x_{n})
\end{equation*}
where $0<\Xi_j<+\infty$ is the normalizing constant, $\phi$ non-negative and measurable. Then 
\begin{equation*}
\Poi=\overset{\infty}{\underset{j=1}{\bigast}} \Qp_j.
\end{equation*}
}

\subsection{Properties of the branching Pólya sum process \label{sect:properties}}

Up to this point we were concerned with the disintegration of the Campbell measure of the Pólya processes with respect to the process itself, i.e. the discussion of the Papangelou kernels. The disintegration with respect to the intensity measure yields further insight. In this section we shall see that it is convenient to add the index $\rho$ to the branching Pólya sum process $\Poi$ to emphasize the dependence on $\rho$. The parameter $z$ stays fixed in this section in contrast to the following section about the Martin-Dynkin boundary.

For the Pólya sum process itself the Palm kernel was identified in~\cite{mR11jtp} as
\equ{
  \PoyS[x]{z,\rho}=\Poy_{z,\rho}\ast\frac{1-z}{z}\sum_{j\geq 1}z^j\Delta_x^{\ast j}\qquad \nu^1_{\Poy_{z,\rho}}\text{-a.e. }x
}
i.e. the Palm kernel of the Pólya sum process at some $x\in X$ is given by the original process with an additional point with geometrically distributed multiplicity. Rewritten this is
\equ{
  \PoyS[x]{z,\rho}=\Poy_{z,\rho}\ast\Poy_{z,\delta_x}\ast\Delta_x
}
and this way to think about the Palm kernel serves as a way to identify the Palm kernels of the branching Pólya sum process.

Note that intuitively, since $\Poy_{z,\rho}$ realizes multiple points, each of which is transformed according to $\kappa$ independently of the others, the observation of a point of $\Pp_\rho$ at some site $x$ in some region $B\in\Esig$ has an influence on the occurence of points at different sites of $B$. 

\satzn{Palm kernels of the branching Pólya sum process}{
Let $z\in(0,1)$, $\rho\in\MX$ and $\kappa$ be a kernel from $X$ to $X$. Then under ${\mathbf A}_1'$ and ${\mathbf A}_2$ for $\nu^1_{\Poi_\rho}$-a.e. $x\in X$
\equ{
  \Poi_\rho^x=\Poi_\rho\ast\Poi_{\delta_x}\ast\Delta_x.
}
}
\bem{
Indeed, the Palm kernel at $x$ of the branching Pólya process is the original process with an additional point at $x$ and another geometric number of points distributed independently according to $\kappa_x$. Conidering the examples in section~\ref{sect:branchmech}, we have
\begin{enumerate}
  \item If $\Esig=\B$, then $\kappa_x=\delta_x$ and the geometric number of points accumulate at $x$, thus generating the multiplicities. Indeed, in this case $\Poi_{\delta_x}=\Poy_{z,\delta_x}$.
  \item If $X=\R^d$ and $\Esig$ is the $\sigma$-algebra of permutation invariant events, then $\Poi_\rho^x$ is the original process with an additional point at $x$ and a geometric number of points distributed on the set of points whose coordinates are permutations of the coordinates of $x$.
  \item If $X=\R^d$ and $\Esig$ is the $\sigma$-algebra of rotation invariant events, then again $\Poi_\rho^x$ is the original process with an additional point at $x$ and a geometric number of points now distributed on the sphere which is centered at the origin and contains $x$.
  \item If $\Esig$ is generated by a locally finite partition of $X$, then $x$ is contained in exacly one of the $X_j$'s and exactly there the geometric number of points is realized.
\end{enumerate}
}

By using ${\mathbf A}_2$ we identify its intensity measure as
\equ{
  \nu^1_{\Poi_\rho}(g)=C_{\Poi_\rho}(g\otimes 1)=\frac{z}{1-z}\int g(x)\kappa\rho(\d x).
}
What remains to do is to disintegrate the Campbell measure with respect to the intensity measure.

\begin{proof}
A successive application of the partial integration and ${\mathbf A}_1'$ to $\kappa$ yields
\eqa{
  C_{\Poi_\rho}(h) &= \iint \sum_{j\geq 1}z^j h(x_j,\mu+\delta_{x_j}+\ldots+\delta_{x_1}) \kappa_{x_{j-1}}(\d x_j)\cdots\kappa_{x}(\d x_1)\rho(\d x)\Poi_\rho(\d\mu)\\
	&= \begin{multlined}[t]\iint \sum_{j\geq 1}z^j h(x_j,\mu+\delta_{x_j}+\ldots+\delta_{x_1})\\ \times\kappa_{x_j}(\d x_{j-1})\cdots\kappa_{x_j}(\d x_1)\Poi_\rho(\d\mu)\kappa_{x}(\d x_j)\rho(\d x).
		\end{multlined}
}
Now observe that the last integration yields $\kappa\rho$, up to a constant the intensity measure $\nu^1_{\Poi_\rho}$, hence by relabelation we get
\equ{
  C_{\Poi_\rho}(h)= \iint \sum_{j\geq 1}z^j h(x_j,\mu+\delta_{y}+\nu) V_{(j-1)\delta_y}(\d\nu)\Poi_\rho(\d\mu)\kappa\rho(\d y)
}
and identify the branching of $\Poy_{z,\delta_y}$.
\end{proof}

\section{The conditioned branching Pólya sum process \label{sect:limits}}

\subsection{Local specifications and H-sufficient statistics \label{sect:limits:constr}}

At first we recall the fundamental objects as in~\cite{eD78,hF75}, then construct the local specification of our interest and finally give the results.

$\Bbd$ is a partially ordered set when equipped with $\subseteq$. If $X$ can be exhausted by an increasing sequence $(B_n)_{n\in\N}$ of bounded sets such that each $B\in\Bbd$ is contained in one of the $B_n$'s, then the limits, which are at a later point to be specified, do not depend on the particular choice of the sequence of sets exhausting $X$. Given a decreasing family $\Ef$ of $\sigma$-fields $\Esig_B\subseteq\Fsig$ indexed by the bounded sets, a family of Markovian kernels $\gamma=(\gamma_{B})_{B\in\Bbd}$ from $\MpmX$ to $\MpmX$ is called a local specification if
\begin{enumerate}
  \item if $B'\subseteq B$, then $\gamma_{B}\gamma_{B'}=\gamma_{B}$; \label{enum:results:locspec:konsistenz}
  \item if $f$ is $\F$-measurable, then $\gamma_{B}(f)$ is $\Esig_B$-measurable; \label{enum:results:locspec:messbarkeit}
  \item if $f$ is $\Esig_B$-measurable, then $\gamma_{B}(f)=f$. \label{enum:results:locspec:ident}
\end{enumerate}
A point processes $\Pr$ with the property 
\begin{equation} \label{eq:results:invariant}
  \Pr(\phi|\Esig_B)=\gamma_{B}(\,\cdot\,,\phi)\qquad\Pr\fsm
\end{equation}
for all $B\in\Bbd$ is called \emph{stochastic field with local specification $\pi$}. The set $C(\gamma)$ of these stochastic fields is convex and the aim is to derive integral representations of these stochastic fields in terms of the extremal points $\Delta\subseteq C(\gamma)$. Since $B\subseteq B_{n_0}$ for each $B\in\Bbd$ and some $n_0$ depending on $B$, $C(\gamma)$ agrees with the set of all point processes $\Pr$ such that equation~\eqref{eq:results:invariant} holds for each $B_n$.

A $\sigma$-field $\Esig$ is \emph{sufficient} for  $C(\gamma)$ if the $\Pr\in C(\gamma)$ have a common conditional distribution given $\Esig$, i.e. there exists $Q_\mu$ such that
\begin{equation*}
  Q_\mu=\Pr\bigl(\,\cdot\,|\Esig\bigr)(\mu)\qquad \Pr\text{-a.e. }\mu.
\end{equation*}
According to~\cite{eD78}, the tail-$\sigma$-field $\Etail=\bigcap_n\Esig_n$ is a H-sufficient statistic for $C(\gamma)$, i.e. it is a sufficient statistic and $\Pr\bigl(\mu:Q_\mu\in C(\gamma)\bigr)=1$ for all $\Pr\in C(\gamma)$. 

Define the $\sigma$-algebra $\F_B$ of events inside a set $B\in\B$ as
\equ{
  \F_B=\sigma(\zeta_{B'}: B'\in\B, B'\subseteq B),
}
and $\hat{\F}_B=\F_{B^c}$ the events outside $B$. Choose an increasing sequence of invariant events $(B_n)_{n\in\N}$ according to assumption ${\mathbf A}_2$ and let
\equ{
  \F_n=\hat{\F}_{B_n}\vee\sigma(\zeta_{B_n}),
}
then $(\F_n)_n$ is a decreasing sequence of $\sigma$-algebras. Denote by $\Poi_z$ the branching of the Pólya sum or Pólya difference process for a given kernel $\kappa$ with parameters $z$ and $\rho$, the latter one is fixed in this section and therefore omitted. We define a familiy $\gamma=(\gamma_n)_n$ of stochastic kernels by
\equ{
  \gamma_n(\mu,\phi)=\Poi_z\bigl(\phi|\F_n\bigr)(\mu).
}
Indeed, $\gamma$ is a local specification and therefore $\Ftail=\bigcap_{n\in\N}\F_n$ is H-sufficient for the $C(\gamma)$.


\subsection{Results \label{sect:limits:results}}

The only point processes which are in some sense canonical branching Pólya processes are mixtures of branching Pólya processes, the only variable parameter turns out to be this $z$, which is constant in case of trivial mixtures. Moreover, this parameter depends on the global particle density. 

For questions concerning convergence we need to assume a continuity condition on $\kappa$ with regard to a measure $\rho$.
\begin{enumerate}
  \item[${\mathbf A}_3$] $x\mapsto\kappa_x(f)$ is $\rho$-a.s. continuous for bounded and continuous $f$ with bounded support
\end{enumerate}

\satzn{Martin-Dynkin boundary of the branching Pólya sum process}{
Assume that $X$ satisfies ${\mathbf A}_2$ and let $\rho\in\MX$ be a diffuse and infinite measure. If $\kappa$ be a kernel from $X$ to $X$ satisfying ${\mathbf A}_1'$ with respect to the given $\sigma$-algebra $\Esig$ and $\rho$ and $\kappa$ satisfy ${\mathbf A}_3$, then the tail-$\sigma$-algebra $\F_\infty$ is a H-sufficient statistic for the family
\equ{
  C(\gamma)=\left\{ \iint V_\tau\,\Poy_{z,\rho}(\d\tau)v(\d z): 
  							v \text{ probability measure on } [0,1) \right\}
}
and its set of extremal points is the family of branching Pólya sum processes indexed by $z\in[0,1)$.
}
Analogously,
\satzn{Martin-Dynkin boundary of the branching Pólya difference process}{
Assume that $X$ satisfies ${\mathbf A}_2$ and let $\rho\in\MX$ be an infinite point measure. If $\kappa$ be a kernel from $X$ to $X$ satisfying ${\mathbf A}_1'$ with respect to the given $\sigma$-algebra $\Esig$ and $\rho$ and $\kappa$ satisfy ${\mathbf A}_3$, then the tail-$\sigma$-algebra $\F_\infty$ is a H-sufficient statistic for the family
\equ{
  C(\gamma)=\left\{ \iint V_\tau\,\GP_{z,\rho}(\d\tau)v(\d z): 
  							v \text{ probability measure on } [0,+\infty) \right\}
}
and its set of extremal points is the family of branching Pólya difference processes indexed by $z\in[0,+\infty)$.
}
\bem{
\begin{enumerate}
  \item Note that one does not need independent increments for all choices of disjoint sets in $B$, as $\Poy_{z,\rho}$ and $\GP_{z,\rho}$. One only needs the independence of the increments for all choices of disjoint sets in $\Esig$. In particular, each branching of a Pólya process has this $\Esig$-independence property and one may use themselves as directing point processes.
  \item Since $\Poy_{z,\rho}$ is infinitely divisible, it can be regarded as a Poisson cluster process with clustering $\Pi_x=\tfrac{z}{1-z}\sum_{j\geq 1}z^j\Delta_x^{\ast j}$. Therefore also the branching Pólya sum process is a Poisson cluster process. However, the family of branching Pólya sum processes is special in the sense that we have a sufficient statistic and knowing some infinite configuration of points, we can estimate the parameter $z$ perfectly.
  \item Usually the continuity of $x\mapsto\kappa_x(f)$ is assumed. Unfortunately this assumption excludes examples of type~~\ref{bsp:setup}~\ref{bsp:setup:part}, where $\kappa(f)$ is not continuous for continuous $f$ in general.
\end{enumerate}
}

The main part of the proofs rely on the identification on the limits $Q_\mu$. Let
\equ{
  U_n=\frac{1}{\rho(B_n)}\zeta_{B_n}
}
be the density of particles in $B_n$ and $U$ be its limit as $n\to\infty$ if it exists. Then
\prop{ \label{thm:limits:results:sum}
Let $X$ satisfy assumption ${\mathbf A}_2$, $\rho$ be an infinite and diffuse measure on $X$, $\kappa$ be a kernel from $X$ to $X$ satisfying ${\mathbf A}_1'$ with respect to the $\sigma$-algebra $\Esig$ and assume that $\rho$ and $\kappa$ satisfy ${\mathbf A}_3$. Then for every $\Pr\in C(\gamma)$ the limit of $(U_n)_n$ exists for $\Pr$-a.e. $\mu$ and for every $\F$ measurable, non-negative $\phi$,
\equ{
  \Pr\bigl(\phi|\F_\infty\bigr)(\mu) = Q_\mu(\phi)
    		=\int V_{\tau}(\phi)\Poy_{Z(\mu),\rho}(\d\tau)\qquad \Pr\text{-a.e. }\mu,
}
where $Z$ is a $\F_\infty$-measurable random variable given by
\equ{
  Z=\frac{U}{1+U}.
}
}
In the same way we get the intermediate result for the branching difference process. Note that in this case $U<1$.
\prop{ \label{thm:limits:results:diff}
Let $X$ satisfy assumption ${\mathbf A}_2$, $\rho$ be an infinite point measure on $X$ and $\kappa$ be a kernel from $X$ to $X$ satisfying ${\mathbf A}_1'$ with respect to the $\sigma$-algebra $\Esig$. Furthermore assume that $\rho$ and $\kappa$ satisfy ${\mathbf A}_3$. Then for every $\Pr\in C(\gamma)$ the limit of $(U_n)_n$ exists for $\Pr$-a.e. $\mu$ and for every $\F$ measurable, non-negative $\phi$,
\equ{
  \Pr\bigl(\phi|\F_\infty\bigr)(\mu) = Q_\mu(\phi)
    		=\int V_{\tau}(\phi)\GP_{Z(\mu),\rho}(\d\tau)\qquad \Pr\text{-a.e. }\mu,
}
where $Z$ is a $\F_\infty$-measurable random variable given by
\equ{
  Z=\frac{U}{1-U}.
}
}
\bem{
In both cases, the sequence $(U_n)_n$ satisfies a strong law of large numbers under each $\Pr\in C(\gamma)$. $U$ is constant a.s. if and only if $\Pr$ is an extremal point of $C(\pi)$.
}

For a probability measure $v$ on $[0,1)$ or $[0,+\infty)$ denote by $\Poi_v$ the mixture of the branching Pólya processes with respect to $v$.

\korollar{ \label{thm:dbstpsp:papangelou}
Let $v$ be a distribution on $(0,1)$ and $\rho_0$ be some infinite and diffuse measure on $\MX$. Then $\Poi_v$ is a solution of the partial integration formula
\equn{ \label{eq:dbstpsp:papangelou}
  C_P(h)=\iint h(x,\mu+\delta_x) Z(\mu)\bigl(\kappa\rho_0+\kappa\mu\bigr)(\d x)P(\d\mu)
}
and $Z$ is the $\Fsig_\infty$-measurable random variable as in Proposition~\ref{thm:limits:results:sum}. Moreover, any solution of~\eqref{eq:dbstpsp:papangelou} with $\Fsig_\infty$-measurable $Z\in(0,1)$ and infinite and diffuse measure $\rho_0\in\MX$ is a mixed branching Pólya sum process and therefore contained in $C(\gamma)$.
}

\korollar{ \label{thm:dbstpdp:papangelou}
Let $v$ be a distribution on $(0,+\infty)$ for some infinite point measure $\rho_0\in\MpmX$. Then $\Poi_v$ is a solution of the partial integration formula
\equn{ \label{eq:dbstpdp:papangelou}
  C_P(h)=\iint h(x,\mu+\delta_x) Z(\mu)\bigl(\kappa\rho_0-\kappa\mu\bigr)(\d x)P(\d\mu)
}
and $Z$ is $\Fsig_\infty$-measurable random variable. Moreover, any solution of~\eqref{eq:dbstpdp:papangelou} with $\Fsig_\infty$-measurable $Z\in(0,+\infty)$ and infinite point measure $\rho_0\in\MX$ is a mixed branching Pólya difference process.
}


\section{Proofs of the construction of the branching Pólya sum process \label{sect:constrproofs}}

Firstly we prove the equivalence of $\pi_\kappa$ satisfying the cocycle condition and $\kappa$ being a regular version of a conditional probability.
\begin{proof}[Proof of Lemma~\ref{thm:equiv-cocycle}]
A direct computation shows that if $\pi$ is either the Pólya sum or the Pólya difference kernel, then for non-negative, measurable functions $f_1$, $f_2$,
\begin{equation*}
  \pi_{\kappa}^{(2)}(\mu;f_{1}\otimes f_{2})= z^{2} \left( \prod_{i=1}^{2} \Bigl(\rho\pm\mu\Bigr)\bigl(\kappa(f_{i})\bigr)\pm\Bigl(\rho\pm\mu\Bigr)\bigl(\kappa(f_{1}\kappa(f_{2}))\bigr) \right).
\end{equation*}
Therefore $\pi_{\kappa}^{(2)}(\mu;f_{1}\otimes f_{2})=\pi_{\kappa}^{(2)}(\mu;f_{2}\otimes f_{1})$ for all $\mu\in\MpmX$ is equivalent to $\kappa(f_{1} \kappa(f_{2}))=\kappa(f_{2} \kappa(f_{1}))$.

Assume firstly that $\kappa$ is a regular version of a conditional probability. Then certainly $\kappa(f \kappa(g))=\kappa(f) \kappa(g)$ holds which implies the cocycle condition.

Assume now that $\pi$ satisfies the cocylce condition. Let us recall Bahadurs \cite{rB55} characterization of conditional expectation: Let $H$ be a probability on $(X,\mathcal{B})$ and $T$ an operator on $L^{2}(X,\mathcal{B},H)$ into itself. Then $T$ is a conditional expectation if and only if:
\begin{displaymath}
 \begin{array}{ll}
 (i) \,\, T \text{ is linear} \,\, (ii)\,\, f\geq 0 \Rightarrow Tf\geq 0 \,\,(iii)\,\, T^{2}=T \\
(iv)\,\, H\left(f\, Tg \right)=H\left(g\, Tf \right) \,\, (v)\,\, T \text{ is constant preserving}.
 \end{array}
\end{displaymath}
Fix $x_{0}\in X$ arbitrarily and define $H=\kappa_{x_{0}}$. We have to check that $f\mapsto \kappa(f)$ satisfies the above properties. From $\kappa(f_{1} \kappa(f_{2}))=\kappa(f_{2} \kappa(f_{1}))$ we directly obtain $\kappa^{2}=\kappa$. Now given $\kappa_{x_{0}}(f^{2})<\infty$ this implies that also with Jensens inequality $\kappa_{x_{0}}(\kappa(f)^{2})\leq \kappa_{x_{0}}(\kappa(f^{2}))=\kappa_{x_{0}}(f^{2})<\infty$. So that $\kappa:L^{2}(X,\mathcal{B},H)\rightarrow L^{2}(X,\mathcal{B},H)$. We certainly have $\kappa_{x_{0}}(f_{1} \kappa(f_{2}))=\kappa_{x_{0}}(f_{2} \kappa(f_{1}))$, which is $(iv)$. $(i),\, (ii),\, (v)$ are clear since $\kappa$ is a stochastic kernel.
\end{proof}

The first step to show that the branching processes are Papangelou processes is to show that they are superpositions of finite processes. Before
\lemma{ \label{thm:brPapkern}
Let $n\in\N$, $B_{1},\ldots, B_{n}\in\Bbd$, $\pi$ be a Papangelou kernel and $\kappa$ satisfying ${\mathbf A}_1$. Then
\equ{
  \pi_{\kappa}^{(n)}(0;B_{1}\times\ldots\times B_{n})
	=\pi^{(n)}(0;\kappa(B_{1})\otimes\ldots\otimes \kappa(B_{n})).
}
}
\begin{proof}
The assertion will be shown by induction. The case $n=1$ holds by definition, and moreover
\begin{align*}
  \lefteqn{\pi_{\kappa}^{(n)}(0;B_{1}\times\ldots \times B_{n})}\qquad\\
	&=\begin{multlined}[t]
	 \iint 1_{B_1}(x_1)\cdots1_{B_{n-1}}(x_{n-1})\pi_{\kappa}(\delta_{x_{1}}+\ldots+ \delta_{x_{n-1}},B_{n})\\
	 \times\pi_{\kappa}^{(n-1)}(0;\d x_{1},\ldots,\d x_{n-1}) 
	\end{multlined}\\
	&=\begin{multlined}[t]
	 \int 1_{B_1}(x_1)\cdots1_{B_{n-1}}(x_{n-1}) \pi(\delta_{x_{1}}+\ldots+ \delta_{x_{n-1}},\kappa(B_{n}))\\
 	\times\kappa_{y_{1}}(\d x_{1})\cdots \kappa_{y_{n-1}}(\d x_{n-1}) \pi^{(n-1)}(0;\d y_{1}\ldots \d y_{n-1})
	\end{multlined}\\
	&=\begin{multlined}[t]
	 \iint 1_{B_1}(x_{1})\cdots 1_{B_{n-1}}(x_{n-1}) \, \kappa_{y_{1}}(\d x_{1})\cdots \kappa_{y_{n-1}}(\d x_{n-1}) \\
	\times\pi(\delta_{y_{1}}+\ldots+ \delta_{y_{n-1}},\kappa(B_n)) \pi^{(n-1)}(0;\d y_{1},\ldots,\d y_{n-1})
	\end{multlined}\\
	&=\pi^{(n)}(0;\kappa(B_1)\otimes\ldots\otimes \kappa(B_n)).
\end{align*}
The third equation follows since for any bounded $B$ and $\mu\in\MpmX$ we have $x\mapsto \pi(\delta_{x}+\mu,\kappa(B))$ is $\Esig$-measurable and $\kappa$ is a regular conditional probability given $\Esig$.
\end{proof}

\begin{proof}[Proof of Proposition~\ref{thm:superpos}]
It suffices to show that the Laplace transforms of both point processes coincide. Denote by $\Poi$ the branching of a Pólya process $\Pp$ and by $X_{(n)}$ the union of $X_1,\ldots ,X_n$. Let $f:X\to[0,\infty)$ measurable such that $\supp f\subseteq X_{(N)}$ for some $N\in \N$. Then since $X_{(N)}\in\Esig$,
\begin{align*}
 \mathcal{L}_{\Poi}(f) &= \int \prod_{x\in \mu} \kappa_{x}(\e^{-f})\Pp(\d\mu)
	= \int \prod_{x\in \mu} \left( \kappa_{x}(1_{X_{(N)}^{c}})+\kappa_{x}(1_{X_{(N)}} \e^{-f}) \right)\Pp(\d\mu)\\
	&=\int \prod_{x\in \mu} \left( 1_{X_{(N)}^{c}}(x)+1_{X_{(N)}}(x) \kappa_{x}( \e^{-f}) \right) \Pp(\d\mu)
	=\int \prod_{x\in \mu_{X_{(N)}}} \kappa_{x}( \e^{-f}) \Pp(\d\mu).
\end{align*}
Now if we denote by $\Pp_{B}$ the restriction of $\Pp$ to $B\in\B$, then $\Pp_{X_{(N)}}=\overset{N}{\underset{j=1}{\bigast}} \Pp_{X_{j}}$, since $\Pp$ has independent increments. Hence we obtain
\equ{
  \mathcal{L}_{\Poi}(f)=\prod_{j=1}^{N} \int \prod_{x\in \nu} \kappa_{x}(\e^{-f})\Pp_{X_{j}}(\d\nu)
  =\prod_{j=1}^{N} \int V_{\nu}(\e^{-\zeta_{f}})\Pp_{X_{j}}(\d\nu).
}
Since as in~\cite{hZ09} locally for $B\in\Bbd$, Pólya sum process and Pólya difference process can be represented as
\equ{
  \Pp_{B}(\phi)=\frac{1}{\Xi(B)} \sum_{n=0}^{\infty} \frac{1}{n!} \int_{B^{n}} \phi(\delta_{x_{1}}+\ldots+\delta_{x_{n}})\, \pi^{(n)}(0;\d x_{1},\ldots,\d x_{n}),
}
where $0< \Xi(B) <\infty$ is a normalizing constant, we get by Lemma~\ref{thm:brPapkern} for $j\in\N$
\equ{
  \Qp_{j}=\int V_{\nu}\,\Pp_{X_{j}}(d\nu).
}
Whence the result follows.
\end{proof}

Since $\pi_{\kappa}$ satisfies the cocycle condition, $\Qp_{j}$ is a finite Papangelou process for each $j$ in the sense of Zessin \cite{hZ09}:
\begin{equation}\label{eq:finite}
  C_{\Qp_{j}}(h)=\int\int_{X_j} h(x,\mu+\delta_{x})\, \pi_{\kappa}(\mu,\d x) \Qp_{j}(\d\mu).
\end{equation}
Now we can identify $\Poi$ as a Papangelou point process with kernel $\pi_{\kappa}$.

\begin{proof}[Proof of Theorem~\ref{thm:papangelou}] 
Let $h$ be of the form $h=f\otimes\e^{-\zeta_{g}}$ with $f,g$ being non-negative, measurable with bounded support, then there exists $N\in \N$ such that $\supp f,\supp g\subset X_{(N)}$, and
\begin{align*}
  C_{\overset{\infty}{\underset{j=1}{\bigast}} \Qp_{j}}(h) 
	&= \int \mu(f) \e^{-\mu(g)}\, \left(\overset{\infty}{\underset{j=1}{\bigast}} \Qp_{j}
 \right)(\d\mu)
	= \int \mu(f) \e^{-\mu(g)}\, \left(\overset{N}{\underset{j=1}{\bigast}} \Qp_{j}\right)(\d\mu)\\
	&= \sum_{j=1}^{N} C_{\Qp_{j}}(h) \prod_{i\in \{1,\ldots,N\}\setminus \{j\}} \mathcal{L}_{\Qp_{i}}(g)\\
	&= \sum_{j=1}^{N} \prod_{i\in \{1,\ldots,N\}\setminus \{j\}} \mathcal{L}_{\Qp_{i}}(g) \int \int_{X_{j}} f(x) \e^{-\zeta_{g}(\nu+\delta_{x})} \, \pi_{\kappa}(\nu_{X_{j}},\d x) \Qp_{j}(\d\nu)\\
	&= \sum_{j=1}^{N} \int \int_{X_{j}} f(x) \e^{-\zeta_{g}(\mu+\delta_{x})} \, \pi_{\kappa}(\mu_{X_{j}},\d x) \left(\overset{N}{\underset{j=1}{\bigast}} \Qp_{j} \right)(\d\mu)\\
	&= \iint f(x) \e^{-\zeta_{g}(\mu+\delta_{x})} \, \pi_{\kappa}(\mu,\d x) 
	 \left(\overset{\infty}{\underset{j=1}{\bigast}} \Qp_{j} \right)(\d\mu).
 \end{align*}
In the fourth equation we used~\eqref{eq:finite}. But the argument in \cite{jM11} chapter 4 proof of Theorem 10 yields that if the integration by parts formula holds for $h=f\otimes \e^{-\zeta_{g}}$, then it holds for every non-negative, measurable $h$. 
\end{proof}

\section{Proofs of the limit theorems for the Pólya sum process \label{sect:limitsproofs}}

In the following we give the proofs of the in section~\ref{sect:limits:results} presented results, the first step is to collect properties of $V_\mu$. Denote by $\Lambda_H$ the unnormalized P\'olya measure for the finite measure $H$ on $X$,
\equ{
  \Lambda_H(\phi)=\sum_{k\geq 0}\frac{1}{k!}\sum_{i_1,\ldots,i_k\geq 1}\frac{z^{i_1+\ldots+i_k}}{i_1\cdots i_k}\int\phi(i_1\delta_{x_1}+\ldots+i_k\delta_{x_n})H(\d x_1)\cdots H(\d x_n).
}
This measure occured in the construction of the P\'olya sum process~\cite{hZ09} as well as in~\cite{mR11jtp}. Note that the unnormalized Poisson measure is defined in a similar way. 
\lemma{ \label{thm:proofs-limit:Vmu}
Let $\Esig\subseteq\B$ be a sub-$\sigma$-algebra and $\kappa_x=H\bigl(\,\cdot\,|\Esig\bigr)(x)$ be a regular version of the probability measure $H$ on $(X,\B)$ conditioned on $\Esig$. Then for every $B\in\Esig\cap\B_0$, $\Lambda_H$-a.e. $\mu\in\MpmX$ and $k\in\N$,
\equ{
  V_\mu(\zeta_{B}=k)=1_{\{\zeta_{B}=k\}}(\mu_{B}).
}
}
For $\Lambda_H$-a.e. point configuration $\mu$ and any invariant set $B$, the number of points of $V_\mu$ realized in $B$ is the same number of point of $\mu$ in $B$. In particular, 
\equ{
  \mu\mapsto V_\mu(\zeta_B=k)
}
is $\sigma(\zeta_B)$-measurable for any $B\in\Esig$.
\begin{proof}
Let $B\in\Esig$, then $\kappa_x(B)=1_B(x)$ $H$-a.e. $x$, that is that a point $x$ in a configuration $\mu$ is contained in $B$, if and only if its image is also contained in $B$. Hence  \eqref{eq:setup:V-def} evaluated at $\{\zeta_{B}=k\}$ for any $k\in\N$ yields
\equ{
  V_\mu(\zeta_{B}=k)=1_{\{\zeta_{B}=k\}}(\mu)\qquad \Lambda_H\text{-a.s},
}
which does not depend on the point configuration $\mu$ outside $B$.
\end{proof}

A direct consequence is that the distribution of $\zeta_B$ under $\Poy_{z,\rho}$ as well as under $\GP_{z,\rho}$ and their branchings agree for $B\in\Esig$.
\korollar{ \label{thm:proofs-limit:zetaB}
Let $\Esig\subseteq\B$ be a sub-$\sigma$-algebra and $\kappa_x=H\bigl(\,\cdot\,|\Esig\bigr)(x)$ be a regular version of the probability measure $H$ on $(X,\B)$ conditioned on $\Esig$ for which $\rho\ll H$. Assume that $X$ satisfies ${\mathbf A}_2'$. Then for every $B\in\Esig$ such that $\rho(B)<+\infty$ and $k\in\N$,
\equ{
  \Poi(\zeta_{B}=k)=\Poy_{z,\rho}(\zeta_{B}=k).
}
The statement remains true if $\Poy_{z,\rho}$ is replaced by $\GP_{z,\rho}$.
}
\begin{proof}
Let $B\in\Esig$ with $\rho(B)<+\infty$, then the restriction of $\Poy_{z,\rho}$ to $B$ is absolutely continuous wrt. $\Lambda_H$ and by Lemma~\ref{thm:proofs-limit:Vmu},
\eqa{
  \iint 1_{\{\zeta_{B}=k\}}(\nu)1_{\{\zeta_{B}=k'\}}(\tau) V_\tau(\d\nu)\Poy_{z,\rho}(\d\tau)
    &=\int 1_{\{\zeta_{B}=k\}}(\tau_{B})1_{\{\zeta_{B}=k'\}}(\tau) \Poy_{z,\rho}(\d\tau)\\
    &=\begin{cases}
    		0 & k\neq k'\\
    		\Poy_{z,\rho}(\zeta_{B}=k) & \text{ otherwise}
    	\end{cases}
}
since if $k=k'$, then the second condition in the first term is superfluous and we get the claim.
\end{proof}

\prop{
Let $\kappa_x=H\bigl(\,\cdot\,|\Esig\bigr)(x)$ be a regular version of the probability measure $H$ on $(X,\B)$ conditioned on the sub-$\sigma$-algebra $\Esig\subseteq\B$, and $\rho\ll H$. Assume ${\mathbf A}_2$, then for every $\F$ measurable, non-negative $\phi$,
\equ{
  \gamma_n(\mu, \phi) = \int \phi(\mu_{B_n^c}+\nu_{B_n})V_{\tau_{B_n}}(\d\nu)
  	      \Pp_{z}\bigl(\d\tau|\F_n\bigr)(\mu),
}
where $\Pp_z$ is either $\Poy_{z,\rho}$ or $\GP_{z,\rho}$.
}
\begin{proof}
Let $A_1\in\hat{\F}_{B_n}$ and $A_2\in\F_{B_n}$ such that $\rho(A_2)>0$. Denote by $\Poi_z$ the P\'olya branching process with parameter $z$. Then because of Corollary~\ref{thm:proofs-limit:zetaB} and the independence of $\F_{B_n}$ and $\hat{\F}_{B_n}$ under $\Pp_{z}$,
\eqa{
  \gamma_n(\mu,A_1\cap A_2) &= 1_{A_1}(\mu)\Poi_z\bigl(A_2|\zeta_{B_n}=\zeta_{B_n}(\mu)\bigr)\\
  		&=1_{A_1}(\mu)\sum_{k\geq 0} \frac{\Poi_z(A_2,\zeta_{B_n}=k)}{\Poi_z(\zeta_{B_n}=k)} 
  					1_{\{\zeta_{B_n}=k\}}(\mu)\\
  		&=1_{A_1}(\mu)\sum_{k\geq 0} \frac{1_{\{\zeta_{B_n}=k\}}(\mu)}{\Pp_z(\zeta_{B_n}=k)} 
  				\iint 1_{A_2}(\nu)1_{\{\zeta_{B_n}=k\}}(\nu)V_\tau(\d\nu)\Pp_{z}(\d\tau)\\
  		&=1_{A_1}(\mu)\sum_{k\geq 0} \frac{1_{\{\zeta_{B_n}=k\}}(\mu)}{\Pp_z(\zeta_{B_n}=k)}
  				\iint 1_{A_2}(\nu)V_{\tau_{B_n}}(\d\nu)1_{\{\zeta_{B_n}=k\}}(\tau) 
  				\Pp_{z}(\d\tau)\\
  		&=1_{A_1}(\mu)\iint 1_{A_2}(\nu)V_{\tau_{B_n}}(\d\nu) 
  				\Pp_{z}\bigl(\d\tau|\zeta_{B_n}=\zeta_{B_n}(\mu)\bigr).
}
Since $\Pp_{z}$ has independent increments, the condition on the Pólya process may be replaced by $\F_n$.
\end{proof}

\begin{proof}[Proof of Proposition~\ref{thm:limits:results:sum}]
Let $f$ be a non-negative, continuous function with bounded support. Then for some $n\in\N$ and then for every larger one, $\supp f\subseteq B_n$, hence $\supp\kappa(f)\subseteq B_n$. Moreover, $\kappa(\e^{-f})$ is bounded away from 0 and therefore $g=-\log\kappa(\e^{-f})$ $\rho$-a.s. bounded and by assumption $\rho$-a.s. continuous. By
\equ{
  \gamma_n(\,\cdot\,,\e^{-\zeta_f}) = \int V_{\tau_{B_n}}(\e^{-\zeta_f})\Poy_{z,\rho}\bigl(\d\nu|\F_n\bigr),
}
since the directing measure converges weakly if and only if the sequence of densities $(U_n)_n$ converges, the limit is $\Poy_{Z,\rho}$ with $Z=\tfrac{U}{1+U}$ by~\cite{mR11jtp}. Therefore we get
\equ{
  Q_\mu(\e^{-\zeta_f}) = \int \e^{-\kappa\nu(f)}\Poy_{Z(\mu),\rho}(\d\nu).
}
\end{proof}

The proof of Proposition~\ref{thm:limits:results:diff} passes the same lines with the help of Theorem~\ref{thm:mdpdp:mdbd} in the next section. For completeness we include the proof of Corollary~\ref{thm:dbstpsp:papangelou}, analogously Corollary~\ref{thm:dbstpdp:papangelou} is proven.

\begin{proof}[Proof of Corollary~\ref{thm:dbstpsp:papangelou}]
Let $v$ be a probability measure on $(0,1)$, then the mixed Pólya sum process $\Poy_v$ solves the partial integration formula~\eqref{eq:dbstpsp:papangelou}, following directly from Proposition~\ref{thm:limits:results:sum}. 

If $P$ is any solution of~\eqref{eq:dbstpsp:papangelou}, then the joint Laplace transform of $Z$ and $P$ for $u,t\geq 0$, $f:X\to\R$ non-negative, bounded and measurable with bounded support is
\equ{
  L_{Z,P}(u,v,tf)=P\left(\e^{-uZ-t\zeta_f}\right)=P\left(\e^{-uZ}P\left(\e^{-t\zeta_f}|\Fsig_\infty\right)\right)
}
by conditioning on $\Fsig_\infty$. Denote by $P_\infty$ the conditioned point process $P$. Differentiation with respect to $t$ yields the Campbell measure of $P_\infty$, which allows to identify this conditional measure, thus on the one hand
\equ{
  -\frac{\d}{\d t}L_{Z,P}(u,tf) = P\left(\e^{-uZ}C_{P_\infty}\left(f\otimes\e^{-t\zeta_f}\right)\right).
}
On the other hand,
\eqa{
  \begin{multlined}
  -\frac{\d}{\d t}L_{Z,P}(u,tf) = C_P\left(f\otimes\e^{-uZ-t\zeta_f}\right)\\
     =P\left(\e^{-uZ}\iint f(x)\e^{-t\mu(f)-tf(x)} Z\bigl(\kappa\rho_0+\kappa\mu\bigr)(\d x)P_\infty(\d\mu)\right)
  \end{multlined}
}
for all $u\geq 0$. Exchanging integration and differentiation is justified since $f\otimes\e^{-uZ-t\zeta_f}$ and $C_{P_\infty}\left(f\otimes\e^{-t\zeta_f}\right)$ are integrable since $f$ is bounded and has bounded support. Thus $P_\infty$ satisfies $P$-a.s. integration by parts formula of the branching Pólya sum process with the parameters given by $Z$ and $\rho_0$ and kernel $\kappa$. But then immediatly $P$ is a mixture of these processes.
\end{proof}

 

\section{Proofs of the limit theorems for the Pólya difference process \label{sect:pdp}}

The identification of the Martin-Dynkin boundary for the Pólya difference process goes along slight different lines than for the Pólya sum process. However, the main partial results are similar. Recall that for an increasing sequence $(B_n)_n$ of bounded sets $(\F_n)_n$ given by
\equ{
  \F_n=\hat{\F}_{B_n}\vee\sigma(\zeta_{B_n}),
}
is a decreasing sequence of $\sigma$-algebras, and the familiy $(\gamma_n)_n$ of stochastic kernels given by
\equ{
  \gamma_n(\mu,\phi)=\GP_{z,\rho}\bigl(\phi|\F_n\bigr)(\mu).
}
for some $z\in(0,+\infty)$ and $\rho\in\MpmX$ is a local specification. Note that $\gamma_n$ does not depend on the parameter $z$. Since $\GP_{z,\rho}$ realizes binomially distributed random numbers at each atom of $\rho$ with the same success probability $\tfrac{z}{1+z}$ by Lemma~\ref{thm:mdpdp:partint} below, there is no surprise that the stochastic fields with that local specification are mixed Pólya difference processes.

\satzn{Martin-Dynkin boundary of the Pólya difference process}{ \label{thm:mdpdp:mdbd}
Let $\rho\in\MpmX$ be an infinite point measure. The tail-$\sigma$-algebra $\F_\infty$ is a H-sufficient statistic for the family
\equ{
  C(\gamma)=\left\{ \int \GP_{z,\rho}\,v(\d z): 
  							v \text{ probability measure on } [0,+\infty) \right\}
}
}

This result follows from the following proposition, which is similar to those in~\cite{mR11jtp}.
\prop{ \label{thm:mdpdp:limits}
Let $\rho\in\MpmX$ be an infinite point measure, $f:X\to\R$ non-negative and measurable with bounded support and $(k_n)_{n}$ be a sequence of poitive integers such that
\equ{
  u=\lim_{n\to\infty}\frac{k_n}{\rho(B_n)}
}
exists with $u\in[0,1)$. Then
\equ{
  \lim_{n\to\infty}\GP_{z,\rho}\left(\e^{-\zeta_f}|\zeta_{B_n}=k_n\right)
  		=\GP_{z',\rho}\left(\e^{-\zeta_f}\right),
}
where $z'=\tfrac{u}{1-u}$.
}

While for the Pólya sum process a representation as a superposition of independent processes was employed, for the Pólya difference process is shown more or less directly. First we give a lemma that adapts the technique in~\cite{hZ09} for the Pólya sum process to the Pólya difference process.

\lemma{ \label{thm:mdpdp:partint}
Let $B\in\Bbd$ and $\phi:\MpmX\to\R$ be non-negative and $\F_{B^c}$-measurable. Then
\equn{ \label{eq:mdpdp:partint}
  \GP_{z,\rho}(\phi)=(1+z)^{\rho(B)}\GP_{z,\rho}(\phi\cdot1_{\{\zeta_B=0\}}).
}
}
\begin{proof}
Let $k\in\N$, $k\leq\rho(B)$, then by applying the integration-by-parts formula and using the $\F_{B^c}$-measurability of $\phi$,
\eqa{
  D_{z,\rho}(\phi\cdot1_{\{\zeta_B=k\}}) 
	&=\frac{1}{k}\iint_B \phi(\mu)1_{\{\mu(B)=k\}}\mu(\d x)D_{z,\rho}(\d\mu)\\
	&=\frac{1}{k}\iint_B \phi(\mu)1_{\{\mu(B)=k-1\}}z\bigl(\rho-\mu\bigr)(\d x)D_{z,\rho}(\d\mu)\\
	&=\frac{z\bigl(\rho(B)-(k-1)\bigr)}{k}\int \phi(\mu)1_{\{\mu(B)=k-1\}}D_{z,\rho}(\d\mu).
}
By iterating one obtains 
\equ{
  \GP_{z,\rho}(\phi\cdot1_{\{\zeta_B=k\}}) 
	=\frac{z^k\rho(B)_{(k)}}{k!}\int \phi(\mu)1_{\{\mu(B)=0\}}\GP_{z,\rho}(\d\mu)
}
for every $k\leq\rho(B)$, and finally by summing equation~\eqref{eq:mdpdp:partint}.
\end{proof}

Directly from this lemma the Laplace transform of the Pólya difference process can be derived.
\korollar{ \label{thm:mdpdp:laplace}
Let $f:X\to\R_+$ be non-negative and measurable with bounded support. Then
\equn{ \label{eq:mdpdp:laplace}
  \GP_{z,\rho}\bigl(\e^{-\zeta_f}\bigr) = \exp\left(\int\log\frac{1+z\e^{-f}}{1+z}\d\rho\right).
}
}
\begin{proof}
Let $B=\supp f$, fix $x\in\supp\rho\cap B$ and decompose $f=f_0+g$ such that $g1_{\{x\}^c},f_01_{\{x\}}=0$ $\rho$-a.s. Then
\eqa{
  \GP_{z,\rho}\bigl(\e^{-\zeta_f}\bigr) &= \GP_{z,\rho}\bigl(\e^{-\zeta_{f_0}}\e^{-\zeta_g}\bigr)\\
	&=\sum_{k\geq 0} \GP_{z,\rho}\bigl(\e^{-\zeta_{f_0}}\e^{-\zeta_g}1_{\{\zeta_{\{x\}}=k\}}\bigr)\\
	&=\sum_{k\geq 0} \e^{-kg(x)}\GP_{z,\rho}\bigl(\e^{-\zeta_{f_0}}1_{\{\zeta_{\{x\}}=k\}}\bigr)\\
	&=\left(1+z\e^{-g(x)}\right)^{\rho(x)}
		\GP_{z,\rho}\bigl(\e^{-\zeta_{f_0}}1_{\{\zeta_{\{x\}}=0\}}\bigr).
}
By iterating one obtains the claimed formula.
\end{proof}

In a similar way we treat the Laplace transform of the conditioned Pólya difference process. In this situation, we choose $n$ such that a given function $f$ with bounded support vanishes outside $B_n$.

\lemma{
Let $f:X\to\R$ be a non-negative and measurable with $\supp f\subseteq B\in\Bbd$. Furthermore let $A\in\Bbd$ with $B\subseteq A$ and choose $m\in\N$ such that $\rho(B)\leq m\leq\rho(A)$. Then
\equ{
  \GP_{z,\rho}\left(\e^{-\zeta_f},\zeta_{A}=m\right) 
    = \GP_{z,\rho}(\zeta_{A}=m)\sum_{k=0}^{\rho(B)} C_{A,m,k}\sum_{\nu\in\Mpm_{(k)}(B):\nu\leq\rho} \prod_{x\in B}{\rho(x)\choose \nu(x)}\e^{-\nu(f)},
}
where
\equ{
  C_{A,m,k}=\frac{m!}{\rho(A)_{(m)}}\cdot\frac{\rho(A\setminus B)_{(m-k)}}{(m-k)!}.
}
}
\begin{proof}
The first step is to distinguish between what happens inside $B$ and what inside $A\setminus B$, and then to apply Lemma~\ref{thm:mdpdp:partint} to $A\setminus B$,
\eqa{
  \GP_{z,\rho}\left(\e^{-\zeta_f},\zeta_A=m\right)
    &=\sum_{k=0}^{\rho(B)} \GP_{z,\rho}\left(\e^{-\zeta_f},\zeta_{A\setminus B}=m-k,\zeta_B=k\right)\\
    &=\sum_{k=0}^{\rho(B)} \frac{z^{m-k}}{(m-k)!} \rho(A\setminus B)_{(m-k)} \GP_{z,\rho}\left(\e^{-\zeta_f},\zeta_{A\setminus B}=0,\zeta_B=k\right).
\intertext{
Next apply the technique used in the proof of Corollary~\ref{thm:mdpdp:laplace} to obtain
}
    &=\begin{multlined}[t]
		\GP_{z,\rho}\left(\zeta_A=0\right)\sum_{k=0}^{\rho(B)} \frac{z^{m-k}}{(m-k)!} \rho(A\setminus B)_{(m-k)}\\
		\times\sum_{k_x\geq0:\sum_x k_x=k}\prod_{x\in B}{\rho(x)\choose k_x}\e^{-k_xf(x)}z^{k_x}
	\end{multlined}\\
    &=\begin{multlined}[t]
		\GP_{z,\rho}\left(\zeta_A=m\right)\sum_{k=0}^{\rho(B)}\frac{\rho(A\setminus B)_{(m-k)}m!}{(n-k)!\rho(A)_{(n)}} \\
		\times\sum_{\nu\in\Mpm_{(k)}(B):\nu\leq\rho}\e^{-\nu(f)}\prod_{x\in B}{\rho(x)\choose \nu(x)}.
	\end{multlined}
}
\end{proof}

\begin{proof}[Proof of Proposition~\ref{thm:mdpdp:limits}]
For sufficiently large $n$, $B_n$ contains $\supp f$. Let $u\in[0,1)$ be the limit of $\tfrac{k_n}{\rho(B_n)}$ as $n\to\infty$ and for simplicity let $A$ be the outer box, $B$ be the inner box and $m$ the number of particles in $A$, then
\eqa{
  C_{A,m,k}&=\frac{m!}{\rho(A)_{(m)}}\cdot\frac{\rho(A\setminus B)_{(m-k)}}{(m-k)!}\\
	&=\frac{m_{(k)}}{\rho(A)_{(k)}}\cdot\frac{\rho(A\setminus B)_{(m-k)}}{\bigl(\rho(A)-k\bigr)_{(m-k)}}.
}
Observe that the first quotient converges to $u^k$, so what remains is the second one. But by applying Stirlings formula, we get asymptotically
\eqa{\allowdisplaybreaks
  \frac{\rho(A\setminus B)_{(m-k)}}{\bigl(\rho(A)-k\bigr)_{(m-k)}}
	&\sim\begin{multlined}[t]
		\sqrt{\frac{\rho(A\setminus B)\bigl(\rho(A)-m\bigr)}{\bigl(\rho(A\setminus B)-m+k\bigr)\bigl(\rho(A)-k\bigr)}}\\
	  \times\frac{\rho(A\setminus B)^{\rho(A\setminus B)}\cdot\bigl(\rho(A)-m\bigr)^{\bigl(\rho(A)-m\bigr)}}{\bigl(\rho(A\setminus B)-m+k\bigr)^{\bigl(\rho(A\setminus B)-m+k\bigr)}\cdot\bigl(\rho(A)-k\bigr)^{\bigl(\rho(A)-k\bigr)}}.\\
	\end{multlined}
}
While the square root tends to 1, the quotient tends to $(1-u)^{\rho(B)-k}$, thus
\equ{
  \lim_{n\to\infty}\GP_{z,\rho}\left(\e^{-\zeta_f}|\zeta_{B_n}\right)
  		=(1-u)^{\rho(B)}\sum_{k=0}^{\rho(B)} \left(\frac{u}{1-u}\right)^k\sum_{\nu\in\Mpm_{(k)}(B):\nu\leq\rho} \prod_{x\in B}{\rho(x)\choose \nu(x)}\e^{-\nu(f)}.
}
Putting $z'=\tfrac{u}{1-u}$, we get that the conditioned Laplace transforms converge to the one given in equation~\eqref{eq:mdpdp:laplace} with $z$ replaced by $z'$.
\end{proof}

\begin{proof}[Proof of Theorem~\ref{thm:mdpdp:mdbd}]
Finally let $\M$ be the set of $\mu\leq\rho$ such that $\lim_n\tfrac{\mu(B_n)}{\rho(B_n)}$ exists as well as $\tilde{\M}$ be the set of $\mu\leq\rho$ such that $\lim_n\gamma_n(\mu,\,\cdot\,)$ exists. Then by Proposition~\ref{thm:mdpdp:limits}
\equ{
  \lim_{n\to\infty}\gamma_n(\mu,\,\cdot\,)=\GP_{Z'(\mu),\rho},
}
hence $\M\subseteq\tilde{\M}$. Vice versa, $\tilde{\M}\subseteq\M$. Therefore we have for any $\Pr\in C(\pi)$,
\equ{
  \Pr(\,\cdot\,|\F_\infty)=\GP_{Z',\rho}\qquad \Pr\text{-a.s}
}
and if $V_\Pr$ is the distribution of $Z'$ under $\Pr$, then
\equ{
  \Pr(\phi)=\Pr\bigl(\GP_{Z',\rho}(\phi)\bigr)=\int \GP_{z',\rho}(\phi)V_\Pr(\d z')
}
and $\Pr$ is a mixed Pólya difference process.
\end{proof}



\bibliographystyle{alpha-abbrv}


\end{document}